\input amstex
\documentstyle{amsppt}
\magnification =1200
\pagewidth {6.5 true in}
\pageheight {8.5 true in }

\parindent 20 pt
\parskip = 3pt plus 1pt minus .5pt
 

\define\fp{\flushpar}
 
\font\bigmath=cmsy10 scaled \magstep1
\def\bigtimes{\mathbin{\hbox{\bigmath\char'2}}}
 
\define\cstar{$\text{C}^*$}  
\define\st{\,|\,}  
\define\BB{\Cal B}
\define\HH{\Cal H}
\define\II{\Cal I}
\define\SS{\Cal S}
\define\orb{\text{orb}}
\define\pred{\operatorname{pred}}
\define\suc{\operatorname{succ}}
\define\ran{\operatorname{ran}}
\define\s{\sigma}
\redefine\t{\tau}
\define\sab{\s_{a,b}}
\define\saa{\s_{a,a}}
\define\tab{\t_{a,b}}
\define\spab{\s_{pa,b}}
\define\ssab{\s_{sa,b}}

\define\phs{\phi_{\s}}
\define\pht{\phi_{\t}}
\define\phab{\phi_{a,b}}
\define\phpab{\phi_{pa,b}}
\define\phasa{\phi_{a,sa}}
\define\pspaab{\psi_{pa,a,b}}
\define\bds#1{(\phs(#1),#1)}
\define\bd#1{(\phi(#1),#1)}
\define\pmax{p_{\text{max}}}
\define\pmin{p_{\text{min}}}
\define\phiminus{\phi^{-}}
\define\phiplus{\phi^{+}}
\define\phat{\phi^{a,t}}
\define\sat{\s^{a,t}}
 
\topmatter
\title
 Boundary Functions for Ideals \\
 in Analytic Limit Algebras
\endtitle
 
\author
Alan Hopenwasser
\endauthor
 
\affil
University of Alabama
\endaffil
 
\address
Department of Mathematics, University of Alabama, Tuscaloosa,
AL 35487
\endaddress
 
\email
ahopenwa\@ua1vm.ua.edu
\endemail
 
\thanks
Much of the research for this paper was done while the
author was on a sabbatical visit to the Technion in
Haifa, Israel.  The author would like to thank the Technion,
and his host, Baruch Solel, for assistance in facilitating
this work.
\endthanks
 
\endtopmatter
\document
\centerline{\bf 1. Introduction}
\par
Boundary functions form a useful tool in the
study of ideals in various classes of nest algebras.  In the simplest
case, where the nest algebra is $T_n$, the algebra of $n \times n$
upper triangular matrices, it is a simple matter to associate to
each ideal in $T_n$ an appropriate boundary function.  This was
generalized to weakly closed ideals in general nest subalgebras of
$\BB(\HH)$ by Erdos and Power in \cite{EP}
 and to Volterra nest subalgebras
of \cstar-algebras by Power in \cite{P1}.
Larson and Solel extended the Erdos-Power theory to the context
of nest subalgebras of factor von Neumann
algebras \cite{LS}.  Both
theories apply to modules over the nest algebra, not just to ideals
in the nest algebra.
Davidson, Donsig and Hudson in \cite{DDH} study support
functions for norm closed bimodules of nest algebras;
their support functions come in pairs which allow the
determination of a maximal and sometimes a minimal
bimodule for a given pair
\par
Amongst algebras which are direct limit of $T_n$'s,
the full nest algebras
introduced in \cite{HP} have the most
in common with weakly closed
nest algebras.  It is not surprising, therefore, that it is
possible to define boundary functions in this context.  The definition
of boundary functions for ideals given in this paper is based on
the possibility of coordinitization for these algebras.  In point
of fact, since everything is based on the properties of the
coordinate system, the theory is valid for a wider class of operator
algebras than full nest algebras.  This class will be the trivially
analytic subalgebras of simple AF \cstar-algebras with an
injective 0-cocycle.
Although the theory of boundary functions in this paper
resembles the theory in the various papers cited in the
first paragraph, one significant difference is that the boundary
functions in the cited papers are all maps from the nest of
invariant subspaces to itself while the boundary functions in
this paper are maps from the spectrum of the diagonal of the
algebra to itself.
\par
We shall assume throughout this paper that $B$ is a simple AF
\cstar-algebra and that $A$ is an analytic subalgebra of $B$
with a trivial cocycle $c$ which is the coboundary of
an injective 0-cocycle.
 These analytic subalgebras are
necessarily maximal triangular. The diagonal,
 $D = A \cap A^*$, of $A$ is a canonical masa
in $B$.
\par
AF \cstar-algebras are groupoid \cstar-algebras; it is this
groupoid and substructures of the groupoid that provide
the necessary coordinitization for the existence of boundary
functions for arbitrary ideals.  We describe these
coordinates very briefly; consult \cite{MS} for
a more detailed treatment.  Since $D$ is an abelian
\cstar-algebra, there is a compact Hausdorff topological space
$X$ such that $D \cong C(X)$.  In the present context,
the spectrum $X$
will, in fact, be a Cantor space.
\par
The groupoid for $B$ will be a principal groupoid based
on $X$; i.e., an equivalence relation on $X$.  One way to
obtain this equivalence relation is as follows: write
$B$ as a direct limit of finite dimensional \cstar-algebras
and choose a system of matrix units for this system.
Each matrix unit in the system acts on $D$ by conjugation
and consequently induces a partial homeomorphism of
$X$, the spectrum of $D$.  The groupoid, $G$, is the
union of the graphs of all the matrix units of the system.
Since the same result is obtained if one uses all
normalizing partial isometries in $B$ instead of a system
of matrix units, $G$ is independent of the choice of
matrix unit system.
\par
In this context, $G$ is an equivalence relation
on $X$ and the groupoid operations are as follows:
\item{(i)}
$(w,x)$ and $(y,z)$ are composible if, and only if, $x=y$,
in which case the product is $(w,z)$, and
\item{(ii)} the inverse
of $(w,x)$ is $(x,w)$.
\flushpar Furthermore, $G$ is a topological
groupoid.  The topology is obtained by declaring that
the graph of each matrix unit will be an open set.  It
turns out that each such graph is also a compact set.
Note that the groupoid
topology is not the relative product topology on
$G$ as a subset of $X \bigtimes X$.
\par
The space, $X$, may be identified with
the diagonal of $X \bigtimes X$ (which
is an open, compact subset of $G$) via
the homeomorphism $x \longrightarrow (x,x)$.
We make this identification hereafter.  Also,
note that the two coordinate projections,
$\pi_1$ and $\pi_2$ of $X \bigtimes X$
onto $X$ are, when restricted to $G$, local
homeomorphisms with respect to the topology on $G$.
An open set on which $\pi_1$ and $\pi_2$ are
one-to-one is known as a $G$-set.  Since any
neighborhood of a point in $G$ contains a smaller
neighborhood which is a $G$-set, we will always
assume in the sequel that neighborhoods are $G$-sets.
\par
The subalgebra, $A$, is a canonical subalgebra, and as
such is generated by the matrix units which it
contains (after a matrix unit system for $B$ has
been selected).  This makes it easy to describe the
support set, $P$, for $A$ (a subset of $G$ whose existence
is guaranteed by the spectral theorem for bimodules
\cite{MS}); $P$ is the union of the graphs of those matrix
units which lie in $A$.  This support set has been
called by a variety of terms in the literature; we
shall refer to it as the spectrum of the algebra $A$
(based on the fact that the relationship between
$A$ and $P$ is strongly analogous to the relationship
between $D$ and $X$).  Similarly, we will sometimes
refer to $G$ as the spectrum of $B$ and, by extension,
to the triple $(X, P, G)$ as the spectral triple for
$(D, A, B)$.  By a theorem of Power \cite{P2},
the spectrum is a complete invariant  for isometric
isomorphism   of  triangular subalgebras of
AF \cstar-algebras.
\par
The spectrum, $P$, satisfies several important
properties in the analytic subalgebra case.  For
example, $P \cap P^{-1} = X$ and $P \cup P^{-1} = G$.
Also, $P$ induces a total order on each equivalence
class from $G$.  If $y$ is an element of $X$, let
$\orb_y$ denote the equivalence class of $x$; i.e.,
$\orb_y = \{x \st (x,y) \in G \}$.  The assumption
that $B$ is a simple \cstar-algebra is equivalent
to the property that each orbit, $\orb_y$, is
dense in $X$ \cite{R, p. 112}.
\par
Examples of spectral triples  of the type under
discussion are provided by the refinement algebras.
Here, $X$ will be the space of all sequences $(x_n)$,
where each term $x_n$ is an element of a set of
positive integers of the form $\{1, \dots , k_n\}$.
Now $X$ is, in fact, the Cartesian product of
countably many finite sets; the topology on $X$
is the product topology for this Cartesian
product.  The equivalence relation, $G$, is the
following: $(x,y) \in G$ if, and only if, there
is an integer $N$ such that $x_n = y_n$ for all
$n \ge N$.  (This is the ``tails are the same''
equivalence relation.)  Each pair of finite
sequences, $(a_1, \dots , a_N)$ and
$(b_1, \dots , b_N)$, determines a basic open set
for the topology in $G$; this set is
$\{ (x,y) \in G \st x_i = a_i \text{ and } y_i = b_i
\text{ for } i = 1, \dots ,N\}$.  Finally, the
spectrum, $P$, is given as follows:
$(x,y) \in P$ if, and only if, there is an integer,
$N$,  such that $(x_1, \dots, x_N) \preceq
(y_1, \dots, y_N)$  (in the lexicographic order) and
$x_i = y_i$ for  $i > N$.
The spectral triple for any refinement algebra can be
represented as described above.
\par
We mention in passing that the spectral triple for any
full nest algebra has a similar representation.  The only
change that is needed is to replace the ``tails the same''
equivalence relation by a possibly much more complicated
equivalence relation which can be determined from a
presentation for the algebra.
\par
Just as the subalgebra, $A$, of $B$ has a support set
$P$ contained in $G$, so does any ideal in $A$ or,
for that matter, any $A$-module, $M$, in $B$.
Suppose that $\II \subseteq A$ is an ideal.  Then
$\II$ is canonical and so is generated by the matrix
units which it contains.  Let $\s$ be the union of
the graphs of the matrix units in $\II$.
The open set $\s$ is the support set of $\II$.
(Since $\II$ is a $D$-bimodule, the spectral theorem
for bimodules may be used to obtain the existence of $\s$.)
The same discussion applies to an $A$-module, $M$.
The support set for an ideal or an $A$-module
satisfies the following definition:
\definition{Definition} An open subset $\s$ of $P$ is
an {\it ideal set\/} if $(w,x), (y,z) \in P$
and $(x,y) \in \s$ imply that $(w,z) \in \s$.
If $\s$ satisfies the same condition but is merely
contained in $G$, then it is an {\it $A$-module set\/}.
\enddefinition
There is a one-to-one correspondence between ideals in $A$
and ideal sets in $P$  and between $A$-modules in $B$
and $A$-module sets in $G$.
Consequently, in the rest of this
paper, we shall discuss ideal sets
(or, once, $A$-module sets) only.
\par
A 1-cocycle, $c$, is a continuous groupoid homomorphism from
$G$ into the real numbers, $\Bbb R$, such that
$c^{-1}(0) = X$.  (Keep in mind that $X$ is
identified with the diagonal of $G$.) The cocycle
property asserts that $c(x,z) = c(x,y) + c(y,z)$
for all $(x,y), (y,z) \in G$.  The canonical
subalgebra $A$ is analytic in $B$ if $P = c^{-1}
[0, \infty)$.  Finally, $A$ is trivially analytic
if $c$ can be written in the form $c(x,y) = b(y) - b(x)$
for some continuous function $b \: X \longrightarrow X$.
As stated above, all the algebras under consideration
are analytic subalgebras of an AF \cstar-algebra which
possess a trivial cocycle.
\par
For example, suppose that $(X,P,G)$ is the spectral
triple of a refinement algebra and has the form
described above.  Then let
$$
\tilde c(x,y) = \sum^{\infty}_{n=1} \frac {y_n - x_n}
  {k_1 k_2 \dots k_n}.
$$
The 1-cocycle, $\tilde c$, is the coboundary of the 0-cocycle
$\tilde b\: X \longrightarrow \Bbb R$ defined by
$$
\tilde b(x) = \sum_{n=1}^{\infty} \frac {x_n -1}
  {k_1 k_2 \dots k_n}.
$$
It is easy to see that $\tilde b$ is a continuous
function on $X$ and that $\tilde c(x,y) = \tilde b(y)
 - \tilde b(x)$, for all $(x,y) \in P$.
\par
We know that $P$ induces a total order on each equivalence
class; in fact, there is a total order on $X$ which agrees
with $P$ on equivalence classes.  This is the lexicographic
order: $x \preceq y$ if, and only if, $x=y$ or there is an
integer $N$ such that $(x_1, \dots, x_N) \prec
(y_1, \dots, y_N)$ in the lexicographic order.
[Note that if $(x_1, \dots, x_N) \prec (y_1, \dots, y_N)$
and $M >N$, then $(x_1, \dots, x_M) \prec
(y_1, \dots, y_M)$.]  The lexicographic order on $X$ has
countably many gaps; let $a^{(n)}$ be an enumeration
of the points with an immediate successor and let
$b^{(n)}$ be the immediate successor of $a^{(n)}$.
These are precisely the points where $\tilde b$ fails
to be one-to-one.  However, we can define a
new and injective 0-cocycle,
$b$, which  induces a 1-cocycle, $c$, given by
$c(x,y) = b(y) - b(x)$.  If $x \in X$,
let $S(x) = \{n \st a^{(n)} \prec x\}$.
Define $b \: X \longrightarrow \Bbb R$ by
$$
b(x) = \tilde b(x) + \sum_{n \in S(x)} \frac 1{2^n}.
$$
Since the order topology induced by $\preceq$ on $X$ is the
same as the topology which $X$ carries as the spectrum
of $D$, $b$ is continuous.  If $(x,y) \in P$, then
$S(x) \subseteq S(y)$ and $c^{-1}[0, \infty) =
 \tilde c^{-1}[0,\infty)$; thus $c$ is a cocycle
for the same analytic subalgebra as $\tilde c$ is.
\par
This same discussion can be carried out in the case
in which $(X,P,G)$ is the spectral triple for a full
nest algebra.  The only change is that the ``equal tails''
equivalence relation is replaced by a more complicated
equivalence relation.  In the full nest algebra case it
is possible that $(a^{(n)}, b^{(n)}) \in P$, for certain
values of $n$.  When this occurs, the formula for $\tilde c$
given above is modified so that $\tilde c(a^{(n)}, b^{(n)}) >0$;
consequently $\tilde b (b^{(n)}) > \tilde b (a^{(n)})$
and the corresponding term may be omitted in the formula
for $b$.
\par
The point of this discussion is that for the class of
algebras of primary interest, the full nest algebras,
there exists an injective 0-cocycle, $b$, on $X$ whose
coboundary, $c$, renders the algebra trivially analytic.
Whenever a trivially analytic algebra is induced from
an injective 0-cocycle, it is possible to define boundary
functions for ideals (or $A$-modules).  Henceforth, we
assume the existence of such an injective 0-cocycle.
\par
The 0-cocycle, $b$, then induces a total order on $X$:
$x \preceq y$ if, and only if $b(x) \le b(y)$.  The
two principal properties of this order are that it
agrees on each equivalence class with the total order
induced by $P$ and that the order topology induced
by $\preceq$ agrees with the original topology
on $X$.  This total order on $X$ is the feature which
makes it possible to define, for each ideal set, a boundary
function $\phi$.
  A simple list of properties characterizes those
functions from $X$ to itself which arise as boundary functions
and there are natural partitions on the family of ideal sets
and on the family of boundary functions so that the quotients
are in bijective  correspondence. (The mapping from
ideal sets to boundary functions is surjective, but
not injective.)
\par
Before proceeding to the definition of boundary functions, we
recapitulate the properties of the spectral
triple which are critical to the notion of a boundary function.
Keep in mind that $G$ has a principal groupoid structure and that
$X$ is identified with the diagonal in $G$, i.e. with the
units of the groupoid.
\item{1.} $G$ is a topological equivalence relation based
on the compact Hausdorff space $X$.
\item{2.} $P$ is an open subset of $G$ which satisfies the
properties: $P \circ P \subseteq P$, $P \cap P^{-1} = X$,
and $P \cup P^{-1} =G$.
\item{3.} The two projection maps, $\pi_1$ and $\pi_2$ of
$X \bigtimes X$ onto $X$ are, when restricted to $G$,
local homeomorphisms with respect to the topology on $G$.
In particular, they are continuous and open mappings.
\item{4.} Each equivalence class, $\orb_y$, from $G$
is countable and dense in $X$.
\item{5.} There is a total order $\preceq$ on $X$ which,
on each equivalence class, $\orb_y$, agrees with the
order induced by $P$.  Furthermore, the order topology
on $X$ is the same as the original topology on $X$.
\item{6.} Since $X$ is compact, it has a minimal
element and a maximal element with respect to
the order $\preceq$.  Denote these elements by
$p_{\text{min}}$ and $p_{\text{max}}$.
\fp
This list of properties will suffice for the definition
and properties of boundary functions.  While there are
spectral triples which satisfy these properties which do
not come from trivially analytic subalgebras of AF
\cstar-algebras with injective 0-cocycle,
the one's that the author knows of are associated with
algebras which lack tractable properties and which do
not appear to be of any interest.  If one adds one further
property -- the assumption that $X$ has at most countably
many gaps -- then it is possible to prove the existence
of a trivial 1-cocycle for $G$ with injective 0-cocycle.
\par
\vskip 12 pt
\centerline{\bf 2. Boundary Functions}
\par
With the preliminaries out of the way, we now turn
attention to the definition and properties of boundary
functions.  Assume that $\s$ is an ideal set contained in
$P$.  Let $y \in X$.  Divide $\orb_y$ into two disjoint
subsets as follows:
$$
\align
A &= \{ x \in \orb_y \st  (x,y) \in \s \}, \\
B &= \{ x \in \orb_y \st (x,y) \notin \s \}.
\endalign
$$
It follows from the definition of ideal set that $A$ is
an initial segment of $\orb_y$ and $B$ is a terminal
segment of $\orb_y$ (in the order induced by $P$).
Now view $A$ and $B$ as subsets of $X$ and
let $v = \sup A$ and $w = \inf B$, where sup
and inf are interpreted with respect to the total
order on $X$ (which agrees with $P$ on $\orb_y)$.
The existence of the sup and inf is guaranteed by
the fact that $X$ is compact in the order topology.
Observe that the open order interval with endpoints
$v$ and $w$ is empty.  Indeed, if this open interval
were not empty, then it would have to contain points
of $\orb_y$ (all equivalence classes are dense in $X$)
and this would contradict the obvious fact that $X$ is
the union of $A$ and $B$.
This leaves two possibilities: either $v=w$ or $v$ is
the immediate predecessor of $w$ in the order, $\preceq$,
on $X$.
\par
In general,
when $x \prec y$ in $X$ and the open order interval with
endpoints $x$ and $y$  is empty, we will say that
$x$ {\it has a gap above\/} and that $y$ {\it has a
gap below\/}.  We will also write $x = \pred y$ and
$y = \suc x$.
\par
\definition{Definition} Let $\s$ be an ideal set in $P$.
 The {\it boundary function for\/} $\s$ is the function
 $\phs \: X \longrightarrow X$ is
given by the formula $\phs(y) = \sup\{x \in \orb_y \st
(x,y) \in \s \}$.
\enddefinition
\par
Note that $\phs(y)$ satisfies the following:
\item{(i)} if $(w,y) \in P$ and $w \prec \phs(y)$, then
 $(w,y) \in \s$, and
\item{(ii)} if $ \phs(y) \prec w$ then $(w,y) \notin \s$.
\fp
It is not possible to say anything about $\bds y$
itself; $\bds y$ may or may not be an element of $P$ and,
if it is an element of $P$, it may or may not be an
element of $\s$.
\par
\proclaim{Proposition 1}
If $\phi = \phs$ is the boundary function for an ideal
$\s$, then $\phi$ has the following properties:
\vskip 1 pt
\item{1.} $\phi(y) \preceq y$, for all $y \in X$.
\item{2.} When $\phi(y)$ has a gap below, the following hold:
\itemitem{a.} $\bd y \in P$.
\itemitem{b.} $y$ has a gap below.
\itemitem{c.} $y \prec z \Longrightarrow \phi(y) \prec \phi(z)$.
\itemitem{d.} there is a neighborhood $N$ (a $G$-set) of $\bd y$
such that $(s,t) \in N \Longrightarrow s \preceq \phi(t)$.
\item{3.} If $y \prec z$ then $\phi(y) \preceq \phi(z)$,
for all $y,z \in X$.
\item{4.} If $y$ does not have a gap below, then
$\phi(y) = \sup \{ \phi(t) \st t \prec y \}$.
\endproclaim
\par
\definition{Definition} A function $\phi \: X \longrightarrow X$
which satisfies properties 1) through 4) in Proposition 1
 will be called a
{\it boundary function\/}.
\enddefinition
\remark{Remark} If $\s$ is an $A$-module set rather than
an ideal set, a boundary function for $\s$ may be defined
in precisely the same way as for ideal sets.  Boundary
functions for $A$-module sets satisfy the properties in
Proposition 1 with two exceptions: condition 1) must be
dropped and condition 2a) must be changed to $\bd y \in G$.
With these two changes, the general notion of a boundary
function can again be defined and the theory described below
remains valid, provided some obvious trivial changes are made.
In the (slightly) modified theory, the sets associated with
boundary functions are, of course, $A$-module sets rather than
ideal sets.  From here on, the exposition will be limited to
ideal sets.
\endremark
\par
\demo{Proof of Proposition 1}
  Property 1) follows immediately from the fact that
$x \preceq y$ for all $(x,y) \in \s$.
\par
Assume that $\phs(y)$ has a gap below.
We must have $\bds y \in \s$, for the other possibility
violates the definition of $\phs$.
In particular, $\bds y \in P$ and 2a) is verified.
\par
 Since $\bds y \in \s$, there is
a neighborhood $N$ of $\bds y$ which is contained in $\s$ and
 is a $G$-set.   We may also assume that $\pi_1(N)$
is contained in the order interval $[\phs (y), \pmax]$
(This order interval is open, since $\phs (y)$ has a gap below;
simply intersect the original  neighborhood with
$[\phs (y), \pmax] \bigtimes X$.)
Now suppose that $y$ does not have a gap below.  Since
$\pi_2$ is an open map, $\pi_2(N)$ contains an open order inteval
whose upper endpoint is $y$.  This implies that there is
a point $(a,b) \in N$ such that $a,b \in \orb_y$ and $b \prec y$.
(This uses, once again, the fact that $\orb_y$ is dense in $X$.)
By the assumptions on $N$, we also have $\phs(y) \prec a$.  The
open order interval with endpoints $\phs (y)$ and $a$ is non-empty
(since $\phs (y)$ does not have a gap above) and therefore contains
a point $z$ from $\orb_y$.  Now observe that
$(z,y) = (z,a) \circ (a,b) \circ (b,y)$ with $(z,a)$ and $(b,y)$
in $P$ and $(a,b)$ in $\s$.  Thus $(z,y) \in \s$.  But $\phs (y)
\prec z$, contradicting the definition of $\phs$.  This proves
that $y$ must have a gap below and condition 2b) is verified.
\par
In order to verify condition 2c), assume that $y \prec z$.
Since both $\phs(y)$ and $y$ have gaps below
and $\bds y \in \s$, there is a
neighborhood  $N$ (a $G$-set) of $\bds y$ which is contained
in $[\phs (y), \pmax] \bigtimes [y, \pmax]$.  In particular,
we can find a point $(a,b) \in N$ such that $a,b \in \orb_z$
and $\phs(y) \prec a$ and $y \prec b \prec z$.  (The
possibility $\phs (y) = a$ is eliminated by the assumption
that $N$ is a $G$-set.)  Thus we have $(a,z) = (a,b) \circ
(b,z)$ with $(a,b) \in \s$ and $(b,z) \in P$.  This shows
that $(a,z) \in \s$ and hence that $a \preceq \phs (z)$.
Since $\phs(y) \prec a$, we have $\phs (y) \prec \phs (z)$,
as desired.
\par
Property 2d) follows from property 2a): since $\bds y \in \s$,
there is a neighborhood $N$ of $\bds y$ which is contained
in $\s$.  If $(s,t) \in N$, then $(s,t) \in \s$ and hence
$s \preceq \phs(t)$.
\par
Assume that property 3) does not hold; i.e., assume that there
are points $y,z \in X$ such that $y \prec z$ and
$\phs(z) \prec \phs (y)$.  By property 2c), $\phs (y)$
does not have a gap below.  Choose an element $a \in \orb_y$
such that $\phs (z) \prec a \prec \phs (y)$.  By the definition
of $\phs$, $(a,y) \in \s$.  So, there is a neighborhood $N$
of $(a,y)$ such that $N \subset \s$ and such that
$(s,t) \in N$ implies $\phs(z) \prec s \prec \phs (y)$
and $t \prec z$.  Choose an element $(s,t) \in N$ with
$s,t \in \orb_z$.  We then have $(s,z) = (s,t) \circ (t,z)$
with $(s,t) \in \s$ and $(t,z) \in P$; hence $(s,z) \in \s$.
But now we have both $\phs(z) \prec s  $ and $s \prec \phs(z)$,
a contradition.  So 3) holds.
\par
Assume that $y$ has no gap below (and hence that $\phs(y)$
also has no gap below).  Property 3) shows that
$t \prec y \Longrightarrow \phs(t) \preceq \phs(y)$;
thus $\sup\{\phs(t) \st t \prec y \} \preceq \phs (y)$.
In order to prove equality, we need to assume that
$w \prec \phs(y)$ and show that there is $t \in X$
such that $t \prec y$ and $w \prec \phs (t)$.
Since $\phs(y)$ has no gap below, there is $x \in \orb_y$
such that $w \prec x \prec \phs(y)$.
By the definition of $\phs$, $(x,y) \in \s$.  Let $N$
be a neighborhood of $(x,y)$ which is contained in $\s$.
Since $y$ has no gap below, there is a point $(s,t) \in N$
such that $t \prec y$ and $w \prec s$.  Since $(s,t) \in \s$,
we have also $s \preceq \phs (t)$; thus $w \prec \phs(t)$
and 4) is verified.  \qed
\enddemo
\par
The mapping $\s \longrightarrow \phs$ from ideal sets
to boundary functions is surjective (as we shall see
later) but is not injective.  To find a simple example
of two ideals with the same boundary function, assume
that $(a,b) \in P$, that $a \ne b$
and that  $a$ does not have a gap
below.  Define
$$
\align
\sab &= \{(x,y) \in P \st \text{\,either }
  x \prec a \text{ or } b \prec y \}, \\
\tab &= \sab \cup \{(a,b)\}.
\endalign
$$
We need to make one other assumption: that $\tab$ is an open
subset of $P$.  In the case of refinement algebras,
this assumption  holds for all points $(a,b) \in P$.
For general full nest algebrs, there may be points
in $P$ for which it fails.  If $\tab$ is
open then it is an ideal set; $\sab$ is always an ideal
set.
\par
Define a function $\psi \:X \longrightarrow X$ as
follows:
$$
\psi(y) =\left \{
\alignedat 2
&y, &\quad &\text{if } y \preceq a, \\
&a, &\quad &\text{if } a \prec y \preceq b, \\
&y, &\quad &\text{if } b \prec y.
\endalignedat \right.
$$
It is a simple matter to check that $\psi$ is the
boundary function for both ideals, $\sab$ and $\tab$.
\par
Next, we consider how to associate ideals to boundary
functions.  So, assume that $\phi \: X \longrightarrow
X$ satisfies the four conditions in the
definition of a boundary function.
We define three subsets of $P$ as follows:
$$
\align
\s(\phi) &= \{(x,y) \in P \st x \prec \phi(y) \}, \\
\eta(\phi) &=\{(x,y) \in P \st x \preceq \phi(y) \}, \\
\s[\phi] &=\{(x,y) \in P \st \text{\,there is a
neighborhood } N \text{ of } (x,y) \text{ with }
N \subseteq \eta(\phi) \}.
\endalign
$$
\par
\proclaim{Proposition 2} The set $\s(\phi)$ is an
ideal set in $P$.
\endproclaim
\demo{Proof} First, we show that $\s(\phi)$ is open.
Let $(x,y) \in \s(\phi)$.  We must find a neighborhood, $N$,
of $(x,y)$ such that $N \subseteq \s(\phi)$.
\par
First, assume that $y$ has a gap below.  The two
order intervals $[\pmin, \phi(y))$ and $[y, \pmax]$
are open subsets of $X$ and $(x,y) \in
[\pmin, \phi(y)) \bigtimes [y, \pmax]$.
Consequently, there is a neighborhood, $N$, of
$(x,y)$ such that $\pi_1(N) \subseteq [\pmin, \phi(y))$
and $\pi_2(N) \subseteq [y, \pmax]$.
If $(w,z) \in N$, we have $w \prec \phi(y)$ and
$y \preceq z$.  By property 3), $\phi(y) \preceq
\phi(z)$.  Thus, $w \prec \phi(z)$ and $(w,z) \in \s(\phi)$.
This shows that $N \subseteq \s(\phi)$.
\par
Now we consider the case when $y$ has no gap below.
Property 2b) implies that $\phi(y)$ has no gap below.
Consequently, there is $s \in X$ such that
$x \prec s \prec \phi(y)$.  By property 4), there
is $t \prec y$ such that $s \prec \phi(t)$.  Now,
the order intervals, $[\pmin, s)$ and $(t, \pmax]$
are open in $X$ and $(x,y) \in [\pmin,s)
\bigtimes (t, \pmax]$.  Consequently, there is
a neighborhood, $N$, of $(x,y)$ such that
$\pi_1(N) \subseteq [\pmin, s)$ and $\pi_2(N)
\subseteq (t, \pmax]$.
Let $(w,z) \in N$.  Then $w \prec s$ and $t \prec z$.
By property 3), $\phi(t) \preceq \phi(z)$.  But
$w \prec s \prec \phi(t)$, so $w \prec \phi(z)$.
Thus $(w,z) \in \s(\phi)$ and $N \subset \s(\phi)$.
This shows that $\s(\phi)$ is an open subset of $P$.
\par
It remains to show that $\s(\phi)$ satisfies the ideal
property.  Assume  $(a,x) \in P$, $(x,y) \in \s(\phi)$, and
$(y,b) \in P$.  Property 3) implies that $\phi(y)
\preceq \phi(b)$.  Thus, we have $a \preceq x \prec
\phi(y) \preceq \phi(b)$; hence $(a,b) \in \s(\phi)$
and $\s(\phi)$ is an ideal set.  \qed
\enddemo
\proclaim{Proposition 3} The set $\s[\phi]$ is an
ideal set in $P$.
\endproclaim
\demo{Proof} Let $(x,y) \in \s[\phi]$.  Then there is a
neighborhood, $N$, of $(x,y)$ such tht $N \subseteq
\eta(\phi)$.  Clearly, any point in $N$ is in
$\s[\phi]$; thus $N \subset \s[\phi]$ and $\s[\phi]$
is an open subset of $P$.
\par
To see that $\s[\phi]$ satisfies the ideal property,
let $(a,x) \in P$, $(x,y) \in \s[\phi]$, and
$(y,b) \in P$.  All neighborhoods in the following argument
are to be open $G$-sets which are subsets of $P$.
Let $N_2$ be a neighborhood of $(x,y)$ such that
$N_2 \subset \eta(\phi)$.  Let $N_1$ be a neighborhood
of $(a,x)$ and $N_3$, a neighborhood of $(y,b)$.
Let $N = N_1 \circ N_2 \circ N_3$.  Then $N$ is
a neighborhood of $(a,b)$.  If $(s,t) \in N$, then
there exist $s', t' \in X$ such that $(s,s') \in N_1$,
$(s',t') \in N_2$, and $(t',t) \in N_3$.  Then,
using $N_1 \subseteq P$, $N_2 \subseteq \eta(\phi)$,
$N_3 \subseteq P$, and property 3) applied to $t' \preceq t$,
we have $s \preceq s' \preceq \phi(t') \preceq \phi(t)$.
Thus $(s,t) \in \eta(\phi)$; since $(s,t)$ is arbitrary
in $N$, $N \subseteq \eta(\phi)$.  This proves that
$(a,b) \in \s[\phi]$, so $\s[\phi]$ satisfies the ideal
property. \qed
\enddemo
\par
We shall see shortly that the boundary function
for the ideal $\s[\phi]$ is $\phi$, thus verifying that
the mapping from ideals to boundary functions is
surjective.  The ideal $\s(\phi)$ need not have $\phi$
as its boundary function.  An examination of some
examples indicates that the boundary function for
$\s(\phi)$ is closely related to $\phi$ and suggests
the following definition.
\par
\definition{Definition} Let $\phi$ be a boundary
function.  Define another function $\phiminus$ by
the formula
$$
\phiminus(y) = \left\{
\alignedat 2
&\phi(y), &\quad &\text{if } \phi(y) \text{ has
no gap below}, \\
&\pred \phi(y), &\quad &\text{if } \phi(y)  \text{ has
a gap below}.
\endalignedat  \right.
$$
\enddefinition
\remark{Remark} When $\phi$ is a boundary function, so
is $\phiminus$.  Property 1) is obvious, since $\phiminus(y)
\preceq \phi(y)$, for all $y$.  There is nothing to prove
for property 2), since $\phiminus(y)$ never has a gap below.
For property 3), assume that $y \prec z$.  If $\phi(z)$ does
not have a gap below, then $\phiminus(y) \preceq \phi(y)
\preceq \phi(z) = \phiminus(z)$.  If $\phi(z)$ does have a
gap below and $\phi(y) \prec \phi(z)$, then $\phiminus(y)
\preceq \phi(y) \preceq \phiminus(z)$.  Finally, if
$\phi(z)$ has a gap below and $\phi(y) = \phi(z)$, then
clearly $\phiminus(y) = \phiminus(z)$.
\par
This leaves 4) to be verified.  When $y$ has no gap below,
neither does $\phi(y)$, by property 2b).  We then have
$$
\align
\phiminus(y) &= \phi(y) = \sup\{\phi(t) \st t \prec y\} \\
&= \sup\{\phi(t) \st t \prec y \text{ and }
 \phi(t) \text{ has
no gap below}\} \\
 &= \sup\{\phiminus(t) \st t \prec y \}.
\endalign
$$
\par
Since $\phiminus(y)$ never has a gap below, we have
$(\phiminus)^{-} = \phiminus$; thus there is never any
need to iterate the ``minus'' operation.
\endremark
\proclaim{Proposition 4} Let $\phi$ be a boundary function.
Let $\s$ be any ideal set such that $\s(\phi) \subseteq
\s \subseteq \s[\phi]$.  Then the boundary function
$\psi$ for $\s$ satisfies $\phiminus \preceq \psi \preceq \phi$.
\endproclaim
\demo{Proof} Let $y \in X$.  First, we show that $\psi(y)
 \preceq \phi(y)$.  We distinguish two cases.  First, assume
that $\psi(y)$ has a gap below.  Then, by property 2a),
$(\psi(y),y) \in \s$.  Since $\s \subseteq \s[\phi]$,
we have
$(\psi(y), y) \in \s[\phi]$ and hence $\psi(y) \preceq
\phi(y)$.
\par
Now assume that $\psi(y)$ has no gap below.  Suppose
that $\phi(y) \prec \psi(y)$.  Since orbits
are dense, there is $t \in \orb_y$
such that $\phi(y) \prec t \prec \psi(y)$.
From the definition of boundary functions for an ideal,
we have $(t,y) \in \s \subseteq \s[\phi]$.
  But this implies that
$t \preceq \phi(y)$, a contradiction.   Thus $\psi(y)
\preceq \phi(y)$ in this case also.
\par
Next we  prove that $\phiminus(y) \preceq \psi(y)$.
Assume, to the contrary, that $\psi(y) \prec \phiminus(y)$.
Since orbits are dense and $\phiminus(y)$ has no gap below,
there is $t \in \orb_y$ such that $\psi(y) \prec t
\prec \phiminus(y)$.  In particular, $t \prec \phi(y)$,
so $(t,y) \in \s(\phi)$.  But $\s(\phi) \subset \s$, so
$(t,y) \in \s$.  The combination $\psi(y) \prec t$ and
$(t,y) \in \s$ contradicts the fact that $\psi$ is the
boundary function for $\s$.  Thus $\phiminus(y)
\preceq \psi(y)$. \qed
\enddemo
\proclaim{Proposition 5} Assume that
  $\s$ is an ideal set and that
$\phi$ is the boundary function for $\s$. Then
$\s(\phi) \subseteq \s \subseteq \s[\phi]$.
\endproclaim
\demo{Proof} Let $(x,y) \in \s(\phi)$.  Then $(x,y) \in P$
and $x \prec \phi(y)$.  From the definition of boundary
function, $(x,y) \in \s$.  Thus $\s(\phi) \subseteq \s$.
\par
Now suppose that $(x,y) \in \s$.  The definition of
boundary function precludes the possibility that
$\phi(y) \prec x$; thus $x \preceq \phi(y)$.  This
shows that $\s \subseteq \eta(\phi)$.  Since $\s$
is open, there is a neighborhood, $N$, of $(x,y)$
such that $N \subseteq \s$.  In particular,
$N \subseteq \eta(\phi)$ and so $(x,y) \in \s[\phi]$.
Thus $\s \subseteq \s[\phi]$.  \qed
\enddemo
The next proposition shows that the mapping from
ideals to the class of boundary functions is surjective.
\proclaim{Proposition 6} Let $\phi$ be a boundary function.
Then the boundary function for the ideal $\s[\phi]$
is $\phi$ and the boundary function for the ideal
$\s(\phi)$ is $\phiminus$.
\endproclaim
\par
\demo{Proof} Let $\psi$ denote the boundary function for
the ideal $\s[\phi]$.  Proposition 4 implies that
$\phiminus \preceq \psi \preceq \phi$.
 If $\phi(y)$ has no gap
below, $\psi(y) = \phi(y)$, so we need only consider the
case in which $\phi(y)$ has a gap below.  Property 2d)
in the definition of boundary functions implies that
$(\phi(y), y) \in \s[\phi]$.  Since $\psi$ is the boundary
function for $\s[\phi]$, we have $\phi(y) \preceq \psi(y)$.
Thus $\psi(y) = \phi(y)$.  This shows that the boundary
function for $\s[\phi]$ is $\phi$.
\par
Now, let $\psi$ denote the boundary function for $\s(\phi)$.
Again, we have that $\phiminus \preceq \psi \preceq \phi$.
Since $\phiminus(y) = \phi(y)$ when $\phi(y)$ has no gap
below, we need only show that $\psi(y) = \phiminus(y)$
whenever $\phi(y)$ has a gap below.
We know from the properties of boundary functions that
$\bd y \in P$.  From the definition of $\s(\phi)$ we
also know that $\bd y \notin \s(\phi)$.  It now follows
that $\psi(y) = \pred \phi(y)$, i.e.~that $\psi(y) =
\phiminus(y)$.  This shows that the boundary function
of $\s(\phi)$ is $\phiminus$. \qed
\enddemo
While the mapping from ideal sets to boundary functions
is not injective, it is possible to say something about
the family of ideals whose boundary function is a
particular function $\phi$.  Of course, all of these
ideals must lie between $\s(\phi)$ and $\s[\phi]$.
The next result says that $\s$ will have $\phi$ for
its boundary function provided that $\s$ contains
an appropriate subset of the graph of $\phi$.
\par
Let $B_{\phi}$ denote the portion of the graph of $\phi$
which is contained in $\s[\phi]$; i.e.,
$$
B_{\phi} = \{\bd y \st y \in X \} \cap \s[\phi].
$$
Let
$$
L_{\phi} = \{\bd y \in B_{\phi} \st \phi(y) \text{ has a gap
below}\;\}.
$$
We then have:
\proclaim{Proposition 7} Let $\phi$ be a boundary function.
An ideal set $\s$ will have $\phi$ for its boundary function
if, and only if, $\s(\phi) \cup L_{\phi} \subseteq \s
\subseteq \s[\phi]$.
\endproclaim
\demo{Proof}  If $\phi$ is the boundary function for $\s$,
then the inclusion $L_{\phi} \subseteq \s$ follows immediately
from the way in which boundary functions are associated with
ideal sets.  (This was pointed out in the proof of part 2a) in
Proposition 1.)  For the converse, assume that
$\s(\phi) \cup L_{\phi} \subseteq
\s \subseteq \s[\phi]$.  Then we know that
 $\phiminus \preceq \phs \preceq \phi$.
We need only check those points $y$ for which $\phiminus(y) \ne
\phi(y)$.  For such $y$, $\bd y \in L_{\phi}$; hence $\bd y
\in \s$.  Therefore, $\phs(y) =\phi(y)$ and $\phs = \phi$ as
desired. \qed
\enddemo
\remark{Remark} Note that the set $\s(\phi) \cup L_{\phi}$
need not be open in  $P$.  In particular, $\s(\phi) \cup
L_{\phi}$ will not, in general, be an ideal set.  As a
consequence, there need not exist a minimal ideal set which
has $\phi$ as its boundary function.  It is not difficult
to produce specific examples of this phenomenon.
\endremark
We obtained $\phiminus$ from $\phi$ basically by replacing
$\phi(y)$ by its immediate predecessor whenever $\phi(y)$
has a gap below.  This suggests defining a function $\phiplus$
in an analogous way, replacing $\phi(y)$ by its immediate
successor whenever $\phi(y)$ has a gap above.
This turns out to be too simplistic, however; doing so will not
produce a boundary function.  For example, if $\phi(y)$
has a gap above and $y$ does not have a gap below, then
any function $\psi$ for which $\psi(y) = \suc \phi(y)$
would fail property 2b) from the definition of boundary
function.  Similar obstacles are presented by properties
2a) and 2d).  The following definition tells just where
we should redefine $\phi(y)$ to obtain $\phiplus$.
\definition{Definition} Let $y \in X$.  Say that $y$ is
a {\it point of modification\/} for $\phi$ if the following hold:
\item{(i)} $y$ has a gap below;
\item{(ii)} $\phi(y)$ has a gap above;
\item{(iii)} there is a neighborhood, $N$, (a $G$-set) of
$(\suc \phi(y), y)$ such that $N \subseteq P$ and, for
all $(s,t) \in N$, \,$s \preceq \phi(t)$ when $\phi(t)$
has no gap above and $s \preceq \suc \phi(t)$ when
$\phi(t)$ does have a gap above.
\enddefinition
\remark{Remark} If $y$ is a point of modification for
$\phi$, then the neighborhood $N$ in condition (iii) may
be selected so that it satisfies
 the additional property that if
$(s,t) \in N$ then $\suc \phi(y) \preceq s$ and $y \preceq t$.
To do so, simply intersect a neighborhood satisfying condition
(iii) with $[\suc \phi(y), \pmax] \bigtimes [y, \pmax]$.
\endremark
\par
We now define a boundary function $\phiplus$ which is larger
than $\phi$ and is closely related to $\phi$.
\definition{Definition} If $\phi$ is a boundary function,
define a function $\phiplus \: X \longrightarrow X$ by
$$
\phiplus(y) = \left\{
\alignedat 2
&\suc \phi(y), &\quad &\text{if } y \text{ is a point
of modification for } \phi, \\
&\phi(y), &\quad &\text{otherwise.}
\endalignedat
\right.
$$
\enddefinition
\proclaim{Lemma 8} Let $\phi$ be a boundary function.
Suppose that $\phi(y) \prec \phiplus(y)$ and $y \prec z$.
Then $\phiplus(y) \prec \phi(z)$.
\endproclaim
\demo{Proof} Choose a neighborhood, $N$, of $(\phiplus(y), y)$
which satisfies both the conditions in the definition of point of
modification and the remark immediately following the definition.
Since $y$ does not have a gap above (it has a gap below),
there is $t \in X$ such that $t \in  \pi_2(N)$,
$y \prec t \prec z$, and $t$ has no gap below.  Let $s$ be
such that $(s,t) \in N$.  Since $t$ is not a point of
modification, we have $\phiplus(y) \preceq s \preceq \phi(t)$.
But $(\phiplus(y), y) \in N$, $y \ne t$, and $N$ is a $G$-set;
so $\phiplus(y) \prec \phi(t)$. Now $t \prec z$ implies
that $\phi(t) \preceq \phi(z)$; hence, $\phiplus(y) \prec
\phi(z)$. \qed
\enddemo
\proclaim{Proposition 9} If $\phi$ is a boundary function
then so is $\phiplus$.
\endproclaim
\demo{Proof} We must show that if $\phi$ satisfies the 4 properties
in the definition of boundary function (see Proposition 1), then so
does $\phiplus$.
\par
Property 1) is automatic, except at points of modification.
If $y$ is a point of modification, then $y$ has a gap below
and $\phi(y)$ has a gap above.  In particular, $\phi(y) \ne y$;
so $\phi(y) \prec y$.  Since $\phiplus(y) = \suc \phi(y)$,
we have $\phiplus(y) \preceq y$.
\par
To verify 2), assume that $\phiplus(y)$ has a gap below.
If $y$ is not a point of modification for $\phi$, then
the four conditions in property 2) hold trivially for
$\phiplus$.  So assume that $y$ is a point of modification.
The first condition, $(\phiplus(y), y) \in P$, follows from
condition (iii) in the definition of point of modification.
Condition 2b), that $y$ has a gap below, is
immediate.  If $y \prec z$ then, by Lemma 8, $\phiplus(y)
\prec \phi(z) \preceq \phiplus(z)$; so condition 2c) holds.
Condition 2d), like condition 2a), follows from property
(iii) in the definition of point of modification.
\par
For the verification of property 3), when $y$ is
not a point of modification,
then $\phiplus(y) = \phi(y) \prec \phi(z) \preceq \phiplus(z)$.
If $y$ is a point of modification, Lemma 8 implies condition
3).
\par
Condition 4) is vacuous at points of modification and trivial
elsewhere, since $t \prec y$ implies $\phi(t) \preceq \phiplus(t)
\preceq \phiplus(y)$.  \qed
\enddemo
\proclaim{Lemma 10} Let $\phi$ be a boundary function.
Then $(\phiplus)^- = \phiminus$ and $(\phiminus)^+ = \phiplus$.
\endproclaim
\demo{Proof} If $\phiplus(y) = \phi(y)$, then $(\phiplus)^-(y)
= \phiminus(y)$ is automatic.  Otherwise, $\phiplus(y)$ is
the immediate successor of $\phi(y)$, in which case both
$(\phiplus)^-(y)$ and $\phiminus(y)$ are equal to $\phi(y)$.
\par
We certainly have $(\phiminus)^+(y) = \phiplus(y)$ when
$\phiminus(y) = \phi(y)$, so assume that $\phiminus(y)$
is unequal to $\phi(y)$ and hence is the immediate predecessor
of $\phi(y)$.  Since $\phiplus(y) = \phi(y)$, we only need
to show that $y$ is a point of modification for $\phiminus$.
We have that $\phiminus(y)$ has a gap above by assumption
and that $y$ has a gap below by property 2b) for the
boundary function $\phi$.  The third condition in the definition
of point of modification  applied to $\phiminus$ follows
from the fact that $\phi$ satisfies condition 2d) in the
definition of boundary function. \qed
\enddemo
\proclaim{Proposition 11} Let $\phi$ and $\psi$ be two
boundary functions.  The the following are equivalent:
\item{A.} $\phiminus = \psi^-$,
\item{B.} $\phiplus = \psi^+$,
\item{C.} $\phiminus \preceq \psi \preceq \phiplus$.
\endproclaim
\demo{Proof} The equivalence of conditions A and B follows
directly from Lemma 10.  If A and B hold, then
$\phiminus = \psi^- \preceq \psi \preceq \psi^+ = \phiplus$,
so condition C holds.  Now, assume that C is valid.
For each $y \in X$, either $\phiminus(y) = \phiplus(y)$
or $\phiminus(y)$ is the immediate predecessor of
$\phiplus(y)$.  Consequently, either $\psi(y) = \phiminus(y)$
or $\psi(y)$ is the immediate successor of $\phiminus(y)$.
In either case, $\psi^-(y) = \phiminus(y)$.  Thus A holds.
\qed
\enddemo
\par
Let $\BB$ denote the family of all boundary functions on $X$.
Define an equivalence relation on $\BB$ as follows:
$\phi \approx \psi$ if, and only if, $\phiminus = \psi^-$.
Proposition 11 implies that the equivalence classes
have the following form: $[\phi] = \{ \psi \st
\phiminus \preceq \psi \preceq \phiplus \}$.
\par
Let $\SS$ denote the family of all ideal sets in $P$.  Define
an equivalence relation on $\SS$: $\s \approx \t$ if, and
only if $\phs^- = \phi_{\tau}^-$.  Equivalence classes can be
identified easily: $[\s] = \{ \t \st \s(\phi_{\s}^-)
\subseteq \t \subseteq \s[\phi_{\s}^+] \}$.
\par
While the mapping $\SS \longrightarrow \BB$ given by
$\s \longrightarrow \phs$ is not surjective, it does
induce a natural bijection of $\SS/\approx$ onto
$\BB/\approx$; viz.~$[\s] \longrightarrow [\phi_{\s}]$.
\par
\vskip 12 pt
\centerline{\bf 3. Examples}
\par
Just prior to Proposition 2 we gave an example of two ideal
sets ($\sab$ and $\tab$) which have the same boundary
function.  We add here a few more very simple examples
which illustrate the properties of boundary functions.
\par
The boundary function for the trivial ideal set,
$\emptyset$, (which corresponds to the trivial
ideal $\II = (0)$) is the function $\phi(y) = \pmin$,
for all $y \in X$.  The boundary function for the
ideal set $\s = P$, (which corresponds to the improper
ideal $\II = A$) is the identity function on $X$.
\par
Let $a \in X$ and let $\s = P \setminus \{(a,a)\}$.
Then $\s$ is an ideal set which corresponds to a
maximal ideal in $A$.  If $a$ has an immediate
predecessor, then the boundary function $\phi$
for $\s$ is given by
$$
\phi(y) = \left\{
\alignedat 2
&\pred a, &\quad &\text{if } y = a, \\
&y, &\quad &\text{otherwise}.
\endalignedat
\right.
$$
If $a$ has no immediate predecessor, then
the boundary function for $\s$ is the identity
function.  Thus, when $a$ has no gap below,
$P$ and $P \setminus \{(a,a)\}$ are another pair
of ideal sets with the same boundary function.
\par
The final example is a variation on the example
preceeding Proposition 2.  For this example we
must assume that $A$ is a refinement algebra
and that $(a,b)$ is a point in $P$ such that
$a$ (and hence $b$) does not have a gap below.
Let
$$
\align
\s' &= \{(x,y) \in \sab \st x \ne y \} \\
\t' &= \{(x,y) \in \tab \st x \ne y \}
\endalign
$$
Then $\s'$ and $\t'$ have the same boundary function,
$\phi$, given by
$$
\phi(y) = \left\{
\alignedat 2
&y, &\quad
 &\text{if } y \text{ has no gap below and either }
 y \prec a \text{ or } b \prec y, \\
&\pred y, &\quad
 &\text{if } y \text{ has a gap below and }
 y \prec a \text{ or } b \prec y, \\
&a, &\quad
 &\text{if } a \preceq y \preceq b.
\endalignedat
\right.
$$
Observe that $\phiminus = \phi$ in this case.  In fact,
we have $\s[\phi] = \s(\phi) \cup \{(a,b)\}$.
\vskip 12 pt
\centerline{\bf 4. Meet and Join Irreducible Boundary Functions}
\par
In \cite{DHHLS} the meet irreducible ideal sets are explicitly
described for algebras with spectral triple $(X, P, G)$ for
which there is a total order on $X$ compatible with $P$.
The description runs as follows: for each pair of
points $a,b \in X$, let
$$
\align
\sab &= \{(x,y) \in P \st x \prec a \text{ or } b \prec y \} \\
\tab &= \sab \cup \{(a,b)\}.
\endalign
$$
While $\sab$ is always an ideal set in $P$, in order for
$\tab$ to be an ideal set we must assume
that $(a,b) \in P$ and that $\tab$ is an open subset
of $P$.  Whenever we use $\tab$, we will assume that
these two conditions are satisfied.
A complete list of all the meet irreducible ideal sets
in P is then given as follows:
\item{1.} $\sab$ if $(a,b)\in P$.
\item{2.} $\sab$ if $(a,b) \notin P$ and there is either
no gap above for $a$ or no gap below for $b$.
\item{3.} $\tab$ if $(a,b) \in P$, there is either no gap
above for $a$ or no gap below for $b$, and $\tab$ is open.
\fp
Later, we will give a description of all the join irreducible
ideal sets.  We shall also see that the boundary function
for an ideal set is meet irreducible or join irreducible (in
an appropriate sense) whenever the ideal set is meet or join
irreducible.
\par
In order to talk about meet and join irreducibility for
boundary functions, we need appropriate lattice operations.
The choices are the obvious ones: $\phi \vee \psi =
\max(\phi, \psi)$ and $\phi \wedge \psi = \min(\phi, \psi)$,
both computed pointwise.
We then have:
\proclaim{Lemma 11}
If $\phi$ and $\psi$ are boundary functions, then so
are $\phi \vee \psi$ and $\phi \wedge \psi$.
\endproclaim
\demo{Proof}
The verification that $\phi \vee \psi$ satisfies the conditions
which define boundary function (given in Proposition 1) is
completely routine.  For $\phi \wedge \psi$, the only conditions
whose verification has some  content are conditions 2c),
2d) and 4). We give arguments for these only.  For convenience,
let $\nu$ denote $\phi \wedge \psi$ and let $y \in X$.  Without
loss of generality, we may assume that $\nu(y) = \phi(y)$.
\par
To verify condition 2c), we assume that $\nu(y)$ has a gap below
and that $z$ is a point in $X$ which satisfies $y \prec z$.
Since $\phi$ is a boundary function we know that $\phi(y)
\prec \phi(z)$; i.e., $\nu(y) \prec \phi(z)$.
From the assumption in the last sentence of the preceeding
paragraph, we also have $\phi(y) \preceq \psi(y)$.  If, in
fact, $\phi(y) = \psi(y)$, then condition 2c) applied to $\psi$
yields $\nu(y) = \psi(y) \prec \psi(z)$.  Thus, in this case,
$\nu(y) \prec \nu(z)$.  So assume that $\phi(y) \prec \psi(y)$.
Use property 3) for $\psi$ to see that
$\nu(y) = \phi(y) \prec \psi(y) \preceq \psi(z)$.  Thus, we
have both $\nu(y) \prec \phi(z)$ and $\nu(y) \prec \psi(z)$,
whence $\nu(y) \prec \nu(z)$.
\par
We continue the assumption that $\nu(y)$ has a gap below
for the verification of condition 2d).  Let $N_1$ be a
neighborhood of $(\phi(y), y)$ such that $(s,t) \in N_1$
impies $s \preceq \phi(t)$.  As before, we have $\phi(y)
\preceq \psi(y)$.  Assume first that $\phi(y) = \psi(y)$.
Then there is a neighborhood $N_2$ of $(\psi(y),y) =
(\phi(y),y)$ such that $(s,t) \in N_2$ implies
$s \preceq \psi(t)$.  Let $N = N_1 \cap N_2$.  Then
$N$ is a neighborhood of $(\nu(y), y)=(\phi(y),y) =(\psi(y),y)$
such that $(s,t) \in N$ implies both $s \preceq \phi(t)$
and $s \preceq \psi(t)$; i.e. $(s,t) \in N$ implies
$s \preceq \nu(t)$.
Now  assume that $\phi(y) \prec \psi(y)$.  With
$N_1$  as above, let $N = N_1 \cap ([\phi(y), \psi(y))
\bigtimes [y, \pmax])$.  Then $N$ is also a neighborhood
of $(\phi(y),y)$.  If $(s,t)\in N$, then $s \prec \psi(y)$
and $y \preceq t$.  By property 3) for $\psi$, we have
$\psi(y) \preceq \psi(t)$; thus $ s \preceq \psi(t)$.
Since $N \subseteq N_1$, we also have $s \preceq \phi(t)$.
This shows that $s \preceq \nu(t)$ and completes the
verification of 2d).
\par
We now turn to the verification of condition 4) for $\nu$
and assume that $y$ does not have a gap below.  For any
$t \prec y$, $\phi(t) \preceq \phi(y)$.  Hence
$\nu(t) \preceq \phi(y) = \nu(y)$ and we have
$\sup\{\nu(y) \st t \prec y\} \preceq \nu(y)$.
Let $x \prec \nu(y) = \phi(y)$.  Since $\phi$ satisfies
condition 4), there is $t_1$ such that
$t_1 \prec y$ and $x \prec \phi(t_1)$.  But $x \prec \psi(y)$
also (since $\phi(y) \preceq \psi(y)$) and $\psi$ satisfies
condition 4); therefore there is $t_2$ such that
$t_2 \prec y$ and $x \prec \psi(t_2)$.  Let
$t_3 = \max(t_1, t_2)$.  Clearly $t_3 \prec y$.
Also, $x \prec \phi(t_1) \preceq \phi(t_3)$ and
$x \prec \psi(t_2) \preceq \psi(t_3)$.  Thus
$x \prec \nu(t_3)$ and $\sup\{ \nu(t) \st t \prec y\}
= \nu(y)$.
\qed
\enddemo
\par
As to be expected, the lattice operations on boundary
functions are related to the lattice operations (set
union and intersection) on ideal sets.  Recall that
the boundary function of an ideal set, $\s$, is
given by $\phs = \sup\{x \in \orb_y \st
(x,y) \in \s\}$.
\proclaim{Lemma 12}
Let $\s$ and $\t$ be ideal sets with boundary functions
$\phs$ and $\pht$.  Then the boundary functions for
the ideal sets $\s \cap \t$ and $\s \cup \t$ are given
by
$$
\align
\phi_{\s \cap \t} &= \phs \wedge \pht \\
\phi_{\s \cup \t} &= \phs \vee \pht
\endalign
$$
\endproclaim
\demo{Proof}
For each $y \in X$ and for each ideal set $\s$,
 $\{x \in \orb_y \st (x,y) \in \s \}$ is an
initial segment of $\orb_y$.
Therefore, for a fixed $y$, the initial segments for $\s$
and for $\t$ are related by inclusion.  Assume, without
loss of generality, that
$$
\{x \in \orb_y \st (x,y) \in \s \} \subseteq
  \{x \in \orb_y \st (x,y) \in \t \};
$$
in other words, assume that $\phs(y) \preceq \pht(y)$.
We then have
$$
\align
\{x \in \orb_y \st (x,y) \in \s \} &=
  \{ x \in \orb_y \st (x,y) \in \s \cap \t \} \\
\{x \in \orb_y \st (x,y) \in \t \} &=
  \{ x \in \orb_y \st (x,y) \in \s \cup \t \}.
\endalign
$$
From this we conclude that $\phs(y) = \phi_{\s \cap \t}(y)$
and $\pht(y) = \phi_{\s \cup \t}(y)$.
But since $\phs(y) \preceq \pht(y)$, we have
$$
\align
(\phs \wedge \pht)(y) &= \phi_{\s \cap \t}(y) \\
(\phs \vee \pht)(y) &= \phi_{\s \cup \t}(y). \qed
\endalign
$$
 \enddemo
In view of Lemma 12, it is natural to expect that the boundary
functions for meet and join irreducible ideal sets are themselves
meet or join irreducible (as appropriate) with respect to
the lattice operations on boundary functions.  First, we consider
the meet operation, for which the following function will
be relevent.  For all $a,b \in X$ with $a \prec b$, define
a function $\phab \:X \longrightarrow X$ by
$$
\phab(y) =  \left\{
\alignedat 2
&y, &\quad &\text{if } y \prec a, \\
&a, &\quad &\text{if } a \preceq y \preceq b, \\
&y, &\quad &\text{if } b \prec y.
\endalignedat  \right.
$$
\par
Provided that $a$ has no gap below, $\phab$ is a boundary
function.  (When $a$ does have a gap below, $\phab$ is not
a boundary function, by property 2b.)  It is straightforward
to check that whenever $\phab$ is a boundary function, it is
a meet irreducible boundary function.
\par
The ideal set $\sab$ is a meet irreducible ideal set except
when $(a,b) \notin P$, $a$ has a gap above, and $b$ has a gap
below.  First, assume that $\sab$ is meet irreducible.
  If $a$ has no gap below, then the boundary function
for $\sab$ is $\phab$.
If $a$ does have a gap below, and if $pa = \pred a$, then
the boundary function for $\sab$ is $\phpab$.  In this case,
$\sab$ and $\spab$ (which
may fail to be a meet irreducible ideal set)
have the same boundary function.
In any event, when $\sab$ is meet irreducible, so is
its boundary function.
\par
If $\sab$ is not meet irreducible, i.e., if
 $(a,b) \notin P$, $a$ has a gap above and $b$ has a gap
below, then $\phab$ (which is meet irreducible) is the boundary
function for $\sab$.  But $\phab$
is also the boundary function of the meet irreducible ideal set
$\ssab$, where $sa = \suc a$.
\par
Thus, whenever $\phab$ is a boundary function (i.e., whenever
$a$ has no gap below), $\phab$ is meet irreducible and the
boundary function of a meet irreducible ideal set.
\par
Next, we consider meet irreducible ideals of the form $\tab$
and their boundary functions.
If $a$ has no gap below, the boundary function for $\tab$
is $\phab$.  (In this case, $\phab$ is the boundary function
of two distinct meet irreducible ideal sets.)  If $a$ has
a gap below, then the boundary function for $\tab$ is
the function $\pspaab$ defined by
$$
\pspaab(y) = \left\{
\alignedat 2
&y, &\quad &\text{if } y \prec a, \\
&pa, &\quad &\text{if } a \preceq y \prec b, \\
&a, &\quad &\text{if } y = b, \\
&y, &\quad &\text{if } b \prec y.
\endalignedat  \right.
$$
It is straightforward to check that $\pspaab$ is
meet irreducible.
\par
Note in passing that if a function of the form $\pspaab$
is a boundary function, then property 2a) implies that
$(a,b) \in P$.  In this case it is also true that $b$
must have a gap below.  (This is required by property
2b) for boundary functions; it is also necessary in order that
$\tab$ be an open set.)
\par
In the case in which $\tab$ is not meet irreducible,
i.e., when $(a,b) \in P$, $a$ has a gap above, and
$b$ has a gap below, the boundary function for $\tab$
is $\phab$.  This function is, of course, meet irreducible
and is also the boundary function of a meet irreducible
ideal set.
\par
In the discussion above, we have assumed that $a \prec b$.
When $a=b$, the ideal set $\saa$ is a maximal ideal set and
hence is meet irreducible.  If $a$ has no gap below, the
boundary function for $\saa$ is the identity function, which
is trivially a meet irreducible boundary function.  If $a$ has
a gap gelow, the boundary function for $\saa$ is $\phasa$, a
meet irreducible boundary function.
\par
\proclaim{Proposition 13} Let $a,b \in P$ with $a \prec b$.
 If $a$ has no gap below, the function $\phab$ defined
above is meet irreducible.  If $a$ has a gap below,
the function $\pspaab$ defined above is meet irreducible.
These functions, together with the identity function, are
the only meet irreducible boundary functions.
Every meet irreducible boundary function is the
boundary function of a meet irreducible ideal set.
Furthermore, if an ideal set is meet irreducible, then
its boundary function is meet irreducible.
\endproclaim
\demo{Proof} We need to show that the boundary functions
listed above are the only meet irreducible boundary
functions.  All the remaining assertions are either
straightforward or have been dealt with in the
discussion preceding the statement of the Proposition.
\par
It is evident from the nature of the meet irreducible
boundary functions, that for a given boundary
function $\phi$, we need to focus on the points
$\phi(y)$ for which $\phi(y) \prec y$.
\define\ED{ED_{\phi}}
\define\rd{RD_{\phi}}
Accordingly, define two sets:
$$
\ED = \{y \st \phi(y) \prec y\} \quad \text{ and } \quad
\rd = \{\phi(y) \st y \in \ED\}.
$$
\par
It is possible that $\rd = \emptyset$.  This happens when
$\phi$ is the identity function.  If $\rd$ is a singleton,
say $\rd = \{a\}$, then $\phi$ has the form $\phab$, for
some $b \in X$.  This is evident from the general fact
that when $a \in \rd$, $\phi^{-1}(a)$ is an order interval
from $X$.  If $\rd$ consists of two points which are the endpoints
of a gap, i.e., if $\rd =\{a,b\}$ where $a$ is the immediate
predecessor of $b$, then $\phi$ has the form $\pspaab$.  In all
of these cases, $\phi$ is a meet irreducible boundary function.
\par
For any other boundary function, $\phi$, there will be two
distinct points in $\rd$ with a third point from $X$ between
the two.  We must show that in this case, $\phi$ is not
meet irreducible.  So assume that $a \prec b \prec c$
in $X$ and that $a,c \in \rd$.
\par
Define an auxiliary function $\eta$ by
$$
\eta(y) = \left\{
\alignedat 2
&y, &\quad &\text{if } y \preceq b, \\
&b, &\quad &\text{if } b \prec y
\endalignedat \right.
$$
and let $\psi_1 = \phi \vee \eta$.
\par
It is evident that $\phi \preceq \psi_1$; furthermore,
$\phi \ne \psi_1$.  Indeed, since $a \in \rd$, there is $z$ such that
$a= \phi(z) \prec z$.  Since $a \prec b$, we have $\phi(z) \prec b$.
Now $\eta(z)$ is either $b$ or $z$.  In either case,
$\phi(z) \prec \eta(z)$.  This means that $\psi_1(z) = \eta(z)
 \ne \phi(z)$.
\par
Now let $t = \sup\{y \in X \st \phi(y) \preceq b \}$.  By properties
3) and 4) for boundary functions, $\phi(t) \preceq b$.  Define
a boundary function $\psi_2$ by
$$
\psi_2(y) = \left\{
\alignedat 2
&\phi(y), &\quad  &\text{if } y \preceq t, \\
&y, &\quad  &\text{if } t \prec y.
\endalignedat  \right.
$$
\par
It is easy to see that $\psi_2$ is a boundary function and
that $\phi \preceq \psi_2$.  Furthermore, $\phi \ne \psi_2$:
there is $s \in X$ such that $c = \phi(s) \prec s$.  Since
$b \prec c = \phi(s)$, $s \notin \{y \st \phi(y) \preceq b\}$.
Since $\phi(t) \preceq b \prec c$, we have $t \prec s$.
Therefore, $\psi_2(s) = s$; in particular, $\psi_2(s) \ne \phi(s)$.
\par
To prove that $\phi$ is not meet irreducible we need
only show that $\phi = \psi_1 \wedge \psi_2$.  Clearly,
$\phi \preceq \psi_1 \wedge \psi_2$.  Let $y \in X$.
If $t \prec y$, then $b \prec \phi(y)$; hence
$\psi_1(y) = \max\{\phi(y), \eta(y)\} = \phi(y)$
(since $\eta(y) \preceq b$).  Thus, $\phi(y) =
(\psi_1 \wedge \psi_2)(y)$.  On the other hand, if
$y \preceq t$, then $\phi(y) = \psi_2(y)$
and $\phi(y) \preceq \psi_1(y)$, so $\phi(y) =
 (\psi_1 \wedge \psi_2)(y)$.  Thus $\phi = \psi_1 \wedge
 \psi_2$. \qed
 \enddemo
\par
Next, we turn to a description of the join irreducible
boundary functions.  Whether or not a boundary function
is join irreducible depends only on the range of the
boundary function.  For any boundary function $\phi$, let
 $\ran \phi  = \{\phi(y) \st y \in X\}$.
Note that, since $\phi(\pmin) = \pmin$, we always
have $\pmin \in \ran \phi$.
\par
\proclaim{Proposition 14}
A boundary function $\phi$ is join irreducible if,
and only if, the cardinality of $\ran \phi$ is at
most 2.
\endproclaim
\demo{Proof}
If $\ran \phi$ contains one element only (necessarily
$\pmin$), then $\phi(y) = \pmin$ for all $y \in X$.
Thus $\phi$ is the minimal boundary function and so
is trivially join irreducible.
\par
Assume that the cardinality of $\ran \phi$ is 2.
Then $\ran \phi = \{\pmin,a\}$, where $\pmin \prec a$.
Observe that $\phi^{-1}(\pmin)$ and $\phi^{-1}(a)$ are intervals
in $X$ which satisfy the property that if $y_1 \in \phi^{-1}(\pmin)$
and $y_2 \in \phi^{-1}(y)$, then $y_1 \prec y_2$.  Furthermore,
the union of these two intervals is all of $X$.  Consequently,
there is an element $t \in X$ such that
$$
\align
y \prec t &\Longrightarrow \phi(y) = \pmin, \\
t \prec y &\Longrightarrow \phi(y) = a.
\endalign
$$
If $t$ does not have a gap below, then it follows from
property 4) of boundary functions that $\phi(t) = \pmin$.
If $t$ does have a gap below, then either alternative,
$\phi(t) = \pmin$ or $\phi(t) =a$ is possible.
Note also that $a \preceq t$, since $\phi(y) \preceq y$, for
all $y$.
\par
Now suppose that $\phi = \psi_1 \vee \psi_2$, where both
$\psi_1$ and $\psi_2$ are boundary function.  It is evident
that on the interval $\phi^{-1}(\pmin)$ we have
$\psi_1 = \psi_2 = \phi$.
\par
First consider the case in which $\phi(t) =a$.  Then
either $\psi_1(t) = a$ or $\psi_2(t) =a$.  Assume,
without loss of generality, that $\psi_1(t) =a$.
Then, for any $y$ with $t \prec y$, we have
$a = \psi_1(t) \preceq \psi_1(y) \preceq a$.  This
shows that $\psi_1(y) =a$ on $\phi^{-1}(a)$, and
thus that $\phi = \psi_1$.
\par
This leaves the case in which $\phi(t) = \pmin$.
Suppose that both $\phi \ne \psi_1$ and
$\phi \ne \psi_2$.  Then there exist elements
$t_1$ and $t_2$ such that $\psi_1(t_1) \prec a$,
$\psi_2(t_2) \prec a$, $t \prec t_1$, and $t \prec t_2$.
Let $t_3 = \min(t_1, t_2)$.  Then $t \prec t_3$,
$\psi_1(t_3) \preceq \psi_1(t_1) \prec a$, and
$\psi_2(t_3) \preceq \psi_2(t_2) \prec a$.  Thus,
$(\psi_1 \vee \psi_2)(t_3) \prec a$ while $\phi(t_3) =a$,
contradicting the assumption that $\phi = \psi_1 \vee \psi_2$.
\par
We have shown that $\phi$ is join irreducible whenever
the cardinality of $\ran \phi$ is at most 2.  We now
assume that the cardinality of $\ran \phi$ is greater
than 2 and show that $\phi$ is not join irreducible.
\par
Assume that $\pmin \prec a \prec b$ and that $a,b \in
\ran \phi$.  We first consider the case in which there
is an element $c \in X$ such that $a \prec c \prec b$.
If there are any points at all between $a$ and $b$, then
there are infinitely many.  In particular, there are points
between $a$ and $b$ with no gap below; so we assume without
loss of generality that $c$ has no gap below.
\par
Let
$$
\align
S &= \{y \st \phi(y) \preceq c \}, \\
T &= \{y \st c \prec \phi(y) \}.
\endalign
$$
Note that $X = S \cup T$ and that $s \in S, t \in T
\Longrightarrow s \prec t$.
\par
Next, define
$$
\align
\psi_1(y) &= \left\{
{\alignedat 2
&\phi(y), &&\quad \text{if } y \in S, \\
&c, &&\quad \text{if } y \in T,
\endalignedat} \right. \\
\psi_2(y) &= \left\{
{\alignedat 2
&\pmin, &&\quad \text{if } y \in S, \\
&\phi(y), &&\quad \text{if } y \in T.
\endalignedat} \right.
\endalign
$$
A routine, but tedious, argument (which we omit) shows
that $\psi_1$ and $\psi_2$ are boundary functions.
Since $b \in \ran \phi$ and $b \notin \ran \psi_1$,
we have $\phi \ne \psi_1$.  Similarly, since
$a$ is in $\ran \phi$ but not in $\ran \psi_2$,
$\phi \ne \psi_2$.  On $S$ it is evident that
$\phi = \psi_1 \vee \psi_2$; since $c \prec \phi(y)$
for all $y \in T$, the same equality is valid on $T$.
Thus $\phi = \psi_1 \vee \psi_2$ and $\phi$ is not
join irreducible.
\par
This leaves the case in which $a,b \in \ran \phi$ and
$b$ is the immediate successor of $a$.  In particular,
$b$ has a gap below.  Let $t$ be such that $\phi(t) =b$.
By condition 2) for boundary functions, $t$ has a gap below.
Furthermore, it $t \prec z$, then $b = \phi(t) \prec \phi(z)$.
So, choose $z$ such that $t \prec z$ (which can be done since
$t \ne \pmax$), and let $d = \phi(z)$.  We now have
$b \prec d$, $b,d \in \ran \phi$ and, since $b$ has no gap above,
there is $c$ such that $b \prec c \prec d$.  By the preceding
argument, $\phi$ is not join irreducible.
\qed \enddemo
\par
\remark{Remark} If $\phi$ is a boundary function whose range
is $\{\pmin,a\}$ with $\pmin \ne a$, then, by property 2),
$a$ cannot have a gap below.  Note that it is also impossible
to have $a = \pmax$.
\endremark
\par
Using Proposition 14, it is easy to describe the join
irreducible boundary functions explicitly.  For each
pair of elements $a,t \in X$ such that $a \preceq t \prec \pmax$,
define $\phat$ by
$$
\phat(y) = \left\{
\alignedat 2
&\pmin, \quad &&\text{if } y \preceq t, \\
&a, \quad &&\text{if } t \prec y.
\endalignedat
\right.
$$
Then $\phat$ is a join irreducible boundary function.
Furthermore, every join irreducible boundary function
is of this form.  (The main issue is the case in which
$t$ has a gap below and $\phi$ is the boundary function
for which $\phi(y) = \pmin$ when $y \prec t$ and $\phi(y)
= a$ when $t \preceq y$.  Let $pt = \pred t$ and note
that, since $a$ has no gap below, $a \prec t$;  in particular,
$a \preceq pt$.  Then $\phi = \phi^{a,pt}$.  The only other
point to note is that $\phi^{\pmin,t}$ is the minimal
boundary function, whose range has cardinality 1.)
\par
If $\phi = \phi^{\pmin,t}$ is the minimal boundary function,
then $\s(\phi) = \s[\phi] = \emptyset$, the ideal set for
the trivial ideal (0). This is the only  ideal set whose
boundary function is the minimal boundary function and it
is trivially a join irreducible ideal set.
\par
For any pair $a,t\in X$ with $\pmin \prec a \preceq t \prec \pmax$,
define an ideal set $\sat$ by
$$
\sat = \{(x,y) \in P \st x \prec a \text{ and } t \prec y \}.
$$
We do not need to assume that $a$ has no gap below; $\sat$
is always an ideal set.  However, the boundary function
for $\sat$ is $\phat$ if, and only if, $a$ has no gap below.
\par
Generally speaking, $\sat$ will be join irreducible.
There is, in fact, only one circumstance when it is not
join irreducible.  This occurs when $a$ has a gap
below (let $pa = \pred a$), $t$ has a gap above
(let $st = \suc t$), and $(pa,st) \notin P$.  In this
case, $\sat = \s^{a,st} \cup \s^{pa,t}$ while
$\sat \ne \s^{a,st}$ and $\sat \ne \s^{pa,t}$.
\par
If $(pa,st) \in P$, then $\sat$ is join irreducible,
as it is in all other cases when either $a$ has no
gap below or $t$ has no gap above.  The verification
that $\sat$ is join irreducible in all these cases
is routine.
\par
If $a$ has a gap below, then the boundary function for
$\sat$ is $\phi^{pa,t}$ (and not $\phat$, which fails
property 2) for boundary functions).  As we shall see
shortly, $\sat$ is the maximal ideal set whose boundary
function is $\phi^{pa,t}$.
\par
Now assume that $a$ has no gap below, so that $\phat$
is a boundary function.  It is evident that $\s(\phat)
= \sat$; thus every join irreducible boundary function
is the boundary function of a join irreducible ideal set.
If $a$ has no gap above, then the properties of boundary
functions ensure that $\s[\phat] = \sat$.  In particular,
when $a$ has no gap above, there is only one ideal set
whose boundary function is $\phat$.  (Use Proposition 7.)
\par
This leaves the case when $a$ does have a gap above.
We then have $\s[\phat] = \s^{sa,t}$.  (Roughly speaking,
because $a$ has a gap above, we can adjoin all the
``boundary points'' $(a,y)$, $y \prec t$ to $\sat$ to
obtain
a set which is open and satisfies the ideal property
 and therefore is an ideal set with the
same boundary function.)
\par
As noted earlier, the only time that $\s[\phat] =
\s^{sa,t}$ will fail to be join irreducible is when
$a$ has a gap above, $t$ has a gap above, and
$(a,st) \notin P$.  Thus, when $a$ has a gap above,
$\phat$ has distinct minimal and maximal ideal sets
amongst the ideal sets whose boundary function is
$\phat$.  The minimal ideal set is always join irreducible
and the maximal ideal set is also join irreducible outside
of one exceptional case.
\par
There are other ideal sets properly between $\s(\phat)$
and $\s[\phat]$ when $a$ has a gap above.  While all of
these have the same join irreducible boundary function,
a routine argument shows that none of these ideal sets
is join irreducible.
\par
This completes the discussion of all ideal sets whose
boundary function has cardinality at most 2.  As for
ideal sets whose boundary function has cardinality
greater than 2, none are join irreducible.  The routine
argument is omitted; it is similar in spirit to the
argument in Proposition 14.
\par
We summarize this discussion as Proposition 15:
\proclaim{Proposition 15} Assume that $a,t \in X$
and $\pmin \prec a \preceq t \prec \pmax$.  Let
$$
\sat = \{(x,y) \in P \st x \prec a \text{ and } t \prec y \}.
$$
Then $\sat$ is a join irreducible ideal set except when
$a$ has a gap below, $t$ has a gap above, and $(pa,st) \notin P$
(where $pa = \pred a$ and $st = \suc t$).  Every non-empty
join irreducible ideal set is of this form.  Every join
irreducible ideal set has a join irreducible boundary function.
Every join irreducible boundary function is the boundary function
of at least one join irreducible ideal set (and at most
two join irreducible ideal sets).
\endproclaim

\enddocument